\documentclass[11pt]{article}

\usepackage[applemac]{inputenc}
\usepackage[T1]{fontenc}

\usepackage{amssymb}
\usepackage[english]{babel}
\usepackage{amsmath}
\usepackage{amscd}

\newcommand{\documentdate}{12 04 2011}
\usepackage{a4wide,graphicx,color}

\let\emptyset\varnothing
\let\geq\geqslant
\let\iff\Leftrightarrow
\let\leq\leqslant
\newcommand{\AND}{\mathbin{\mathrm{and}}}

\newcommand{\Ccs}[1]{\hbox{\mathversion{bold}$\mathrm{Ccs}(#1)$}}
\newcommand{\Cov}{\mathfrak{Z}}
\newbox\bdelta\newbox\bsupscript\newbox\bsubscript
\newcommand{\ddelta}[3]{\setbox\bdelta=\hbox{$\delta$}%
\setbox\bsupscript\hbox{$\scriptstyle #2$}%
\setbox\bsubscript\hbox{$\scriptstyle #3$}
\mathord{\rule{0pt}{\ht\bdelta}%
^{#1\rule[-\dp\bsupscript]{0pt}{1pt}}%
_{\rule{0pt}{\ht\bsubscript}}}\!\delta^{#2}_{#3}}
\newcommand{\eqdef}{\mathrel{:=}}
\def\Fer{\mathop{\mathrm{Fer}}}
\def\HSharp{H^{\sharp}}
\def\id{\mathord{\mathrm{id}}}

\def\Ouv{\mathop{\mathrm{Ouv}}}
\newcommand{\overdot}[1]{\mathord{\stackrel{\circ}{#1}}}
\newcommand{\psdiff}{\stackrel{\mskip6mu\circ}{\leftarrow}}
\newcommand{\Psdiff}{\stackrel{\mskip9mu\circ}{\Leftarrow}}
\newcommand{\PSdiff}[1]{\mathrel{\underset{\mskip10mu\scriptscriptstyle #1}{\stackrel{\mskip10mu\circ}{\Longleftarrow}}}}

\newcommand{\restr}[1]{\mathord{\upharpoonright}_{#1}}

\def\st{\mathbin{|}}
\newcommand{\struct}[1]{\langle #1\rangle}

\def\restoreparskip{\global\parskip 8.4pt plus 1pt \relax} 
\restoreparskip
\belowdisplayshortskip \baselinestretch\baselineskip
\abovedisplayskip\parskip
\belowdisplayskip\parskip
\makeatletter
\def\decalageG{1.8em}
\newcounter{MonEnumeratei}
\newcounter{MonEnumerateii}

\newcount\NiveauItemize \NiveauItemize0

\renewenvironment{itemize}{\advance\NiveauItemize by1\begin{list}%
{\csname EtiquetteItem\romannumeral\the\NiveauItemize\endcsname}%
{\def\makelabel ##1{\hfill##1\hfill}
  \topsep0pt \itemsep0pt \parsep.5\parskip
  \leftmargin\decalageG
  \rightmargin0pt \listparindent0pt \itemindent0pt
  \labelsep0pt \labelwidth\decalageG}%
}{\end{list}}
\renewenvironment{enumerate}{\advance\NiveauItemize by1\begin{list}%
{\csname theMonEnumerate\romannumeral\the\NiveauItemize\endcsname}%
{\usecounter{MonEnumerate\romannumeral\the\NiveauItemize}
  \def\makelabel ##1{\hfill##1\hfill}
  \topsep0pt \itemsep0pt \parsep\parskip \leftmargin\decalageG
  \rightmargin0pt \listparindent0pt \itemindent0pt
  \labelsep0pt \labelwidth\decalageG}%
}{\end{list}}
\renewenvironment{quote}{%
\begin{list}%
{}%
{\topsep0pt \itemsep0pt \parsep\parskip \leftmargin\decalageG
  \rightmargin0pt \listparindent0pt \itemindent0pt
  \labelsep0pt \labelwidth\decalageG}%
  \item[]%
}{\end{list}}
\makeatother
\begin{document}

\begin{titlepage}

\includegraphics[height=3.5cm]{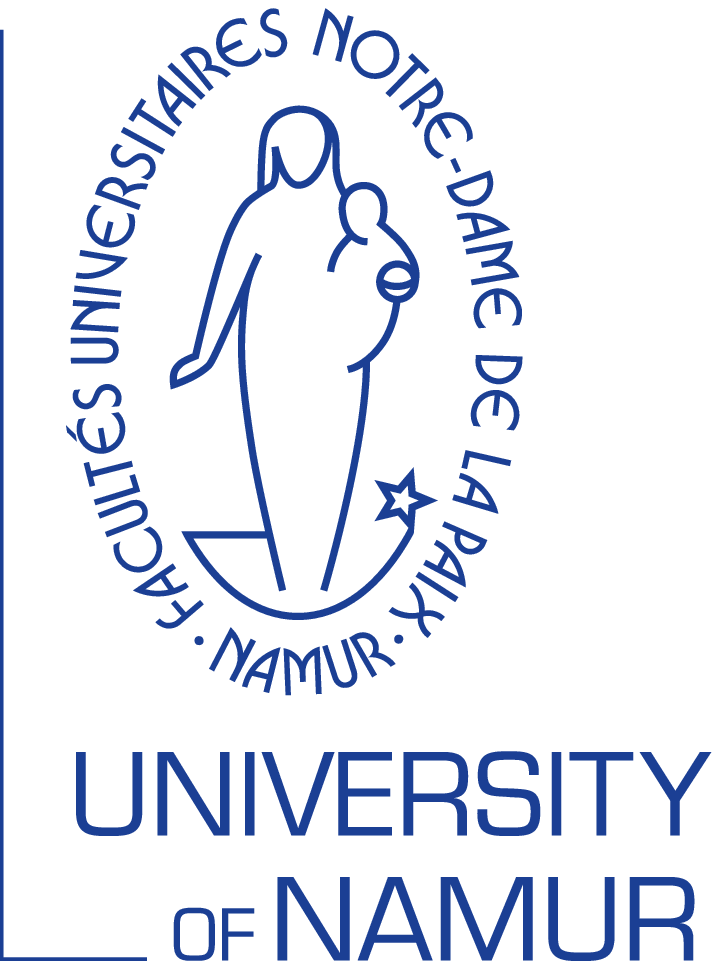}

\vspace*{2cm}
\hspace*{1.3cm}
\fbox{\rule[-3cm]{0cm}{6cm}\begin{minipage}[c]{12cm}
\begin{center}
From contradiction to conciliation:\\
a way to ``dualize'' sheaves \\
\mbox{}\\
by D. Lambert and B. Hespel \\
\mbox{}\\
Report naXys-12-2011 \hspace*{20mm} \documentdate 
\end{center}
\end{minipage}
}

\vspace{2cm}
\begin{center}
\includegraphics[height=3.5cm]{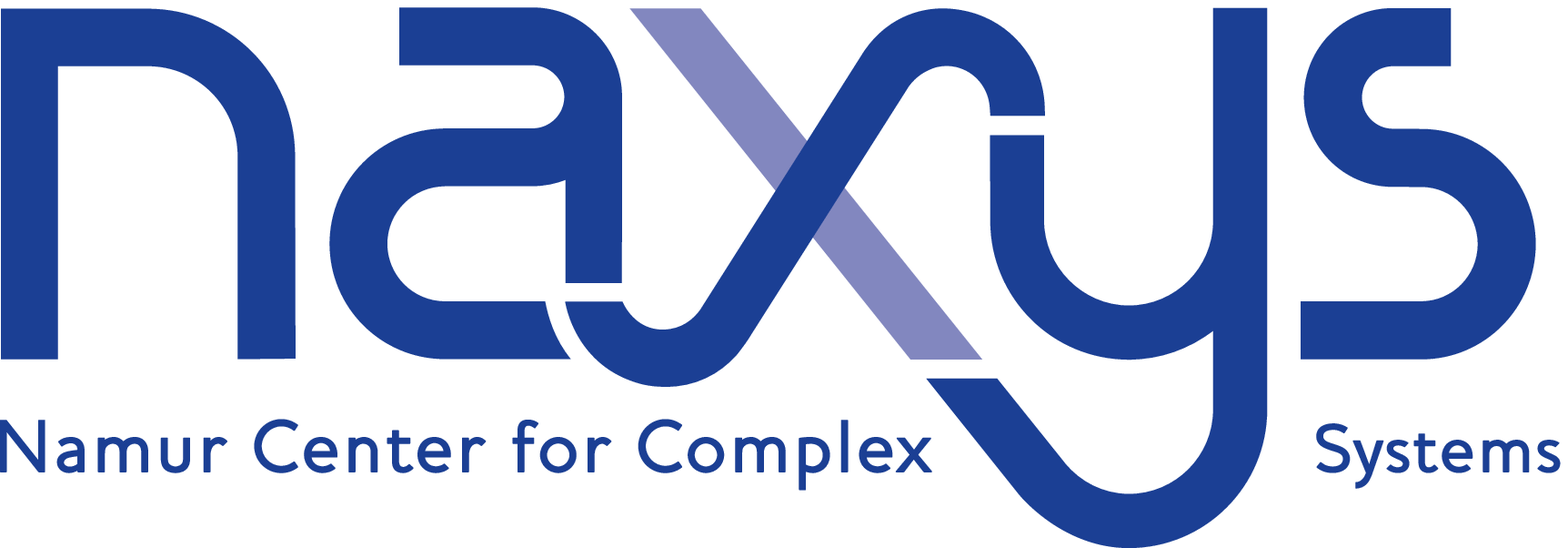}

\vspace{2cm}
{\Large \bf Namur Center for Complex Systems}

{\large
University of Namur\\
8, rempart de la vierge, B5000 Namur (Belgium)\\*[2ex]
{\tt http://www.naxys.be}}

\end{center}

\end{titlepage}
\newpage

\begin{center}
 {\LARGE\bf From contradiction to conciliation:\\[2mm]
 a way to ``dualize'' sheaves \footnote{Both of the authors wish to warmly thank Daniel Dzierzgowski for having greatly improved their manuscript and Professor Chris Mortensen for having quickly shown an interest in their attempt.}}\\[2mm]
 {\large\it D. Lambert and B. Hespel}\\[2mm]
 {\large\ naXys}\\
 {\large\ University of Namur} 
\end{center}

\begin{abstract}
Our aim is to give some insights about how to approach the formal description of situations where one has to conciliate several contradictory statements, rules, laws or ideas. We show that such a conciliation structure can be naturally defined on a topological space endowed with the set of its closed sets and that this specific structure is a kind of “dualization” of the sheaf concept where “global” and “local” levels are reversed.
\end{abstract}

\section{Towards a model of conciliation?}

Our aim here is to give some insights about how to approach the formal description of situations where one has to conciliate several statements, rules, laws, ideas that can be contradictory. 

In fact it is straightforward to notice that many real situations in human social and political life face us with the need to conciliate contradictory points of view and to find some agreement.
Reaching such an agreement does not imply that all contradictions have vanished. The democratic practice shows us that in spite of real contradictions, institutions can survive fruitfully. Of course one realizes easily two things. 
First, the number of contradictory views or positions has to remain under a certain level above which some dangerous tensions could appear with the risk to destroy completely the relations network. 
Thus contradictions have to be confined in some restricted and controlled fields in order to avoid complete percolation in the system under consideration. 
Second, not all contradictory theses can coexist peacefully or without damage for the coherence of institutions, laws or thoughts. 
But, in some circumstances, contradictions survive coherently in a system, institution or in law corpus because they happen to be interpretations, readings, of a common yet ambiguous statement, idea, rule, law,~\ldots\ In fact this means that contradictions can merge in a coherent framework if they are rooted in a common soil from which they can be derived as two non-equivalent but possible interpretations. 

In this paper we try to formalize this kind of realization.
In order to do this, we exhibit a structure called {\it conciliation} taylored for describing situations in which we impose that two different and even contradictory “things” (ideas, statements, rules, images,~\ldots) defined “locally” (i.e. at the level of a restricted social or political groups, law corpus,~\ldots) can be identified (conciliated) at a higher level, “globally” if you want (for example at the level of a country, or of the international law,~\ldots), if and only if there exists at a “deeper level” (for example at the level of citizens, of articles of laws,~\ldots) a unique element which can be map on both initial “things”. 
The philosophy of conciliation building is to impose global constraint (forcing the conciliation of objects at high level) and then to exhibit some objects locally that enable the conciliation. 
The latter are roots providing a deep relation between what can be considered as different or even contradictory at a local level. In fact the existence of such “roots” is not trivial and their existence corresponds precisely to what we mean by a {\it conciliation}.

We will see that such a conciliation structure happens to be naturally defined on a topological space endowed with the set of its closed sets. This conciliation structure will also appear to be a kind of “dualization” of the sheaf concept where “global” and “local” levels are reversed.
In a sheaf, we impose, at a deep level, conditions that enable us to glue coherently objects defined locally in such a way that if the sheaf exists, we can build some new global object whose restrictions give rise to the initial object. 
So it is possible to use this sheaf-conciliation duality to build some new concepts.

\section{The concept of conciliation on a topological space}

Let $X$ be a topological space and $\Fer(X)$ the set of its closed sets.

A {\it pre-conciliation} is defined giving:
\begin{enumerate}
\item a map $G$ which associates to each closed set $W$ belonging
   to $\Fer(X)$ a set $G(W)$;
\item for every $W_1, W_2 \in \Fer(X)$ such that 
   $W_2\supseteq W_1$, a map 
   $\delta^{W_2}_{W_1}: G(W_1) \rightarrow G(W_2)$, called a 
   {\it mediation}, satisfying the following properties:
   \begin{equation}\label{prop-iii}\textstyle
     \forall W\in\Fer(X): \delta^{W}_W=\id_W
   \end{equation}
   where $\id_W$ is the identity map on $W$, and
   \begin{equation}\label{prop-iv}
     \forall W_1, W_2, W_3 \in\Fer(X): 
     W_3\supseteq W_2\supseteq W_1
     \Rightarrow \delta^{W_3}_{W_2} \circ \delta^{W_2}_{W_1}
     = \delta^{W_3}_{W_1}.
   \end{equation}
\end{enumerate}

A family $(W_i)_{i\in I}$ is a {\it finite closed co-covering of 
$W\in \Fer(X)$} iff $I$ is finite, 
$\forall i\in I: W_i\in \Fer(X)$ and 
$W = \bigcap_{i\in I} W_i$.

A pre-conciliation $G$ in $X$ is said to be {\it unified} iff it satisfies
the following property:
\begin{equation}\label{prop-v}
 \begin{array}{l}
  \makebox[11.5cm][l]{for every $W\in \Fer(X)$, $W\neq \emptyset$,}\\[1mm]
  \mbox{and every finite closed co-covering $(W_i)_{i\in I}$ of 
        $W$,}\\[1mm]
  \mbox{we have: $\forall j \in I\;\forall t,s\in G(W):
        \delta^{W_j}_W(t) = \delta^{W_j}_W(s)
        \Rightarrow t=s$.}
 \end{array}
\end{equation}

A unified pre-conciliation $G$ in $X$ is a {\it conciliation} in $X$ iff it satisfies the following property:
\begin{equation}\label{prop-vi}
  \begin{array}{l}
   \makebox[11.5cm][l]{for every $W\in \Fer(X)$, $W\neq \emptyset$,}\\[1mm]
   \mbox{every finite closed co-covering $(W_i)_{i\in I}$ of 
         $W,$}\\[1mm]
  	\mbox{and every $(t_i)_{i\in I}$ such that 
	      $\forall i\in I: t_i\in G(W_i)$,}\\[1mm]
	\mbox{if $\forall i,j \in I: \delta^{W_i\cup W_j}_{W_i}(t_i)
           = \delta^{W_i\cup W_j}_{W_j}(t_j)$}\\[1mm]
    \mbox{then $\exists ! t\in G(W)\;
        \forall i\in I: \delta^{W_i}_W(t)=t_i$.}
  \end{array}
\end{equation}
 
Intuitively, the latter property means the following thing.
The reason why mediations have succeeded to conciliate, at the global level, ideas which are different or even contradictory, is the fact that there exists, at a deep level, a fundamental idea which enables this conciliation. 
And it works because this ground for mutual conciliation, rooting divergent ideas in a common soil.

\section{Some examples of conciliations}

Now we give an example of conciliation that will play a major role to connect this kind of structure to logic.
For any $W_i\in \Fer(X)$, let us consider the map:
\begin{displaymath}
W_i \mapsto \Cov(W_i) = \{Z\in \Fer(X) \st Z \supseteq W_i\}.
\end {displaymath}

Let us define the maps
$\delta^{W_j}_{W_i}:\Cov(W_i)\rightarrow \Cov(W_j):
Z_i\mapsto \delta^{W_j}_{W_i}(Z_i)=Z_i\cup W_j$.
We check that properties (\ref{prop-iii}), (\ref{prop-iv}) and 
(\ref{prop-v}) are trivially satisfied.
We want to show that $\Cov$ is a conciliation in $X$.

Let us suppose that
$\forall i,j\in I \;\forall Z_i\in\Cov(W_i)\;\forall Z_j\in\Cov(W_j):
\delta^{W_i\cup W_j}_{W_i}(Z_i) 
= \delta^{W_i\cup W_j}_{W_j}(Z_j)$,
we have to prove that
$\exists ! Z \in \Cov(W_i\cap W_j):
\delta^{W_i}_{W_i\cap W_j}(Z)=Z_i
\AND \delta^{W_j}_{W_i\cap W_j}(Z)=Z_j$.

By hypothesis, we have:
$\delta^{W_i\cup W_j}_{W_i}(Z_i) 
= \delta^{W_i\cup W_j}_{W_j}(Z_j)
\Leftrightarrow 
Z_i\cup W_i\cup W_j = Z_j\cup W_i\cup W_j$.

But $Z_i\in \Cov(W_i)$ and $Z_j\in\Cov(W_j)$,
therefore we get:
$Z_i\cup W_j = Z_j\cup W_i$.
The searched unique closed set $Z$ is nothing but $Z_i\cap Z_j$.
Indeed, we have
\begin{displaymath}
\begin{array}{c}
    \delta^{W_i}_{W_i\cap W_j}(Z_i\cap Z_j)
    = (Z_i\cap Z_j)\cup W_i
    =(Z_i\cup W_i)\cap (Z_j\cup W_i)
    =Z_i\cap(Z_j\cup W_i)\\
    \delta^{W_j}_{W_i\cap W_j}(Z_i\cap Z_j)
    = (Z_i\cap Z_j)\cup W_j
    =(Z_i\cup W_j)\cap (Z_j\cup W_j)
    =(Z_i\cup W_j)\cap Z_j.   
\end{array}
\end{displaymath}
But we have 
$Z_i\subseteq Z_i \cup W_j = Z_j\cup W_i$,
$Z_j\subseteq Z_j \cup W_i = Z_i\cup W_j$,
then we can deduce that
\begin{displaymath}
 \delta^{W_i}_{W_i\cap W_j}(Z_i\cap Z_j) = Z_i,
 \delta^{W_j}_{W_i\cap W_j}(Z_i\cap Z_j) = Z_j
\end{displaymath}
and, because $\delta^{W_i}_{W_i\cap W_j}$
and $\delta^{W_j}_{W_i\cap W_j}$ are mediations satisfying property
(\ref{prop-v}), that $Z_i\cap Z_j$ is the unique element of 
$\Cov(W_i\cap W_j)$ satisfying these two relations.

A second example of conciliation in $X$ is given considering for any
$W_i\in\Fer(X)$ the map $1:W_i\rightarrow \{*\}\restr{W_i}$ which
associates a singleton to each closed set.
It is sufficient to set trivially:
$\delta^{W_j}_{W_i}(*\restr{W_i}) = *\restr{W_j}$
and we check that this defines a concilation in $X$.

Conciliation has natural links with the concept of extension. We will show this on a third example using the extension theorem of Tietze (see~\cite{Simmons1963}). This states the following:

\begin{quote}
Let $X$ be a normal space and let $F$ one of its closed subsets and $f$ a continuous function defined on $F$ with values in a closed interval
$[a,b]$ of the set of real numbers.\\
Then, $f$ has a continuous extension $g$ defined on $X$ and with values in
$[a,b]$ .
\end{quote}

Let us recall that a topological space is normal if it is a $T_1$ space
in which two disjoint closed sets can always be separated by open sets. 
A $T_1$ space is a topological space in which each point is a closed set. 
We are considering now a compact Hausdorff topological space (a space is a Hausdorff space if two distinct points can always be separated by open sets, in the sense that they have open neighbourhood).
A compact Hausdorff space is a normal space and furthermore every compact subspace of a Hausdorff space is closed.
Let us thus consider two compact subspaces $H_1$ and $H_2$ of $H$ such
that $H_1\cup H_2 = H$.
Now, let $f_1:H_1 \rightarrow [a,b]$ and $f_2:H_2 \rightarrow [a,b]$
be two continuous functions. 
By Tietze extension theorem, we know that both functions admit two continuous extensions: 
$g_1:H \rightarrow [a,b]$ and $g_2:H \rightarrow [a,b]$.
One can write: $f_j = g_j\circ J_j$ ($j=1,2$), where $J_j$ denotes
the inclusion of $H_j$ in $H$.
Let us call $G(H_1)$ and $G(H_2)$ the sets of continuous functions from the closed sets $H_1$ and $H_2$ in $[a,b]$.
We check that extensions of these continuous functions can be viewed as mediations of a conciliation $G$ built on $H_1\cap H_2 = \HSharp$ which
associates to any closed subset a set of continuous functions on this space and with values in $[a,b]$.
We have effectively a conciliation, for if we set the following constraint:
\begin{displaymath}
 g_1=g_2,
 \quad g_1 \eqdef \delta^H_{H_1},
 \quad g_2 \eqdef \delta^H_{H_2},
\end{displaymath}
we check that
\begin{displaymath}
 \exists ! g: H_1\cap H_2 \rightarrow H:
 \delta^{H_1}_{H_1\cap H_2}g = f_1 \;\AND\;
 \delta^{H_2}_{H_1\cap H_2}g = f_2.
\end{displaymath} 
This follows from the fact that
$f_1 \circ K_1 = g_1\circ J_1 \circ K_1
\equiv g_2\circ J_2 \circ K_2 = f_2\circ K_2 \eqdef g$, 
where $K_j$ is the inclusion of $H_1\cap H_2$ in $H_j$ and from the fact that 
$J_1\circ K_1 = J_2 \circ K_2 :H_1\cap H_2\rightarrow H$. 
It is interesting to note that in Tietze theorem, 
the fact that the subspace is closed is crucial. 
If this subspace is not closed, one cannot ensure the existence of a continuous extension and the constraint to get a conciliation is not satisfied (let us take for example 
$H=[0,1]$, $H_1=]0,1]$, 
$f_1:H_1\rightarrow H: x \mapsto f_1(x) = \sin(1/x)$,
the space $H_1$ is not closed in $H$ and we easily check that
$f_1$ has no continuous extension on all of $H$).
This shows that the fact that the definition of conciliation starts from closed sets is not at all artificial. On the contrary the use of closed sets seems dictated by the nature of mediations, whose particular case is extensions.

\section{Conciliations as cosheaves}

Let us note that a pre-conciliation $G$ can be considered as a {\it covariant functor} from the category of closed sets on $X$
(whose objects are the closed sets and whose arrows are closed sets inclusions) in the category of sets.
This functor is described by the following commutative diagram:

\begin{displaymath}
\begin{CD}
W_j @>G>> G(W_j)\\
@A{\textstyle\subseteq}AA @AA{\textstyle\delta^{W_j}_{W_i}}A\\
W_i @>G>> G(W_i)\\
\end{CD}
\end{displaymath}

The above definition of a conciliation in $X$ happens to be a kind of “dualization” of the sheaf concept on $X$ (in that sense a conciliation can be called a “cosheaf”).
We know indeed that a sheaf F on X can be defined as follows (see \cite{MacLane1992}). 

Let $X$ be a topological space and $\Ouv(X)$ the set of its open sets.

A {\it presheaf} $F$ on $X$ is defined giving:
\begin{enumerate}
\item a map $F$ which associates to each open set $U$ of $\Ouv(X)$ a set
  $F(U)$;

\item for every $ U_1,U_2\in\Ouv(X)$ such that $U_1\subseteq U_2$,
  a map $\rho^{W_2}_{W_1}:F(U_2)\rightarrow F(W_1)$,
  called a {\it restriction}, satisfying the following properties:
  \begin{equation}
   \forall U\in\Ouv(X): \rho^{U}_U=\id_U
  \end{equation}
  where $\id_U$ is the identity map on $U$ and
  \begin{equation}
   \forall U_1, U_2, U_3\in \Ouv(X):
     U_1\subseteq U_2\subseteq U_3 \Rightarrow
     \rho^{U_3}_{U_2} \circ \rho^{U_2}_{U_1}
    =\rho^{U_3}_{U_1}.
  \end{equation}
\end{enumerate}

A family $(U_i)_{i\in I}$ is a {\it finite open covering of 
$U\in \Ouv(X)$} iff
$I$ is finite, 
$\forall i\in I: U_i\in \Ouv(X)$ and 
$U = \bigcup_{i\in I} U_i$.

A presheaf on $F$ is said to be {\it separated} if we assume in addition that:
\begin{equation}
 \begin{array}{l}
  \makebox[9.5cm][l]{for every $U\in \Ouv(X))$,}\\[1mm]
  \mbox{and every finite open covering $(U_i)_{i\in I}$
        of $U,$}\\[1mm]
  \mbox{we have: $\forall j \in I\;\forall t,s\in F(U) :
        \rho^U_{U_i}(t) = \rho^U_{U_i}(s)
        \Rightarrow t=s$.}
 \end{array}
\end{equation}

A separated presheaf $F$ on $X$ is a {\it sheaf on $X$} iff it satisfies the following property:
\begin{equation}
  \begin{array}{l}
  \makebox[9.5cm][l]{for every $U\in \Ouv(X))$,}\\[1mm]
  \mbox{every finite open covering $(U_i)_{i\in I}$
        of $U,$}\\[1mm]
  	\mbox{and every $(t_i)_{i\in I}$ such that 
	    $\forall i\in I: t_i\in F(U_i)$,}\\[1mm]
    \mbox{we have: $\exists ! t\in F(U)\,\forall i\in I: 
        \rho^U_{U_i}(t) = t_i$.}
\end{array}
\end{equation}

We know also that a presheaf $F$ on $X$ can be considered as a contravariant functor from the category of open sets (whose objects are open sets and whose arrows are opens sets inclusions) in the category of sets, described by the following commutative diagram:

\begin{displaymath}
 \begin{CD}
  U_j @>F>> F(U_j)\\
  @A{\textstyle\subseteq}AA  @VV{\textstyle\rho^{W_j}_{W_i}}V\\
  U_i @>F>>  F(U_i)
 \end{CD}
\end{displaymath}

We can push forward this kind of “duality” which relates sheaves and conciliations in order to generate new concepts. 

If $F$ is a presheaf on a topological space $X$, we can define at each point
$x$ of $X$ a disjoint union: 
$\bigcup_{U\ni x} F(U)$. 
Now let us define an equivalence relation $\approx_{\rho}$:
for every $U_1, U_2\in \Ouv(X), t_i\in U_i$, 
\begin{displaymath}
 t_1\approx_{\rho} t_2
 \iff 
 \exists U\subseteq U_1\cap U_2: \rho^{U_1}_U(t)
  = \rho^{U_2}_U(t).
\end{displaymath}
In other words, two “sections” of the sheaf are identified if they have the same restriction on an open set. 
Then, classically, one can introduce the concept of stalk 
$F_x$ of $F$ at a point $x$ of $X$ by:
\begin{displaymath}
 F_x \eqdef \bigcup_{U_j\ni x} F(U_j)/\mathord{\approx_{\rho}}
 = \varinjlim_{U_j\ni x}F(U_j)
\end{displaymath}
where the limit is the inductive limit induced by the restrictions: 
\begin{displaymath}
 \begin{array}{l}
  \displaystyle
  F(U_0) \stackrel{\rho^{U_0}_{U_1}}{\longrightarrow} F(U_1)
    \stackrel{\rho^{U_1}_{U_2}}{\longrightarrow} F(U_2)
    \stackrel{\rho^{U_2}_{U_3}}{\longrightarrow} 
    \cdots
    \varinjlim_{U_j\ni x} F(U_j)\\[2mm]
  U_0\supseteq U_1\supseteq U_2\supseteq \ldots
 \end{array}
\end{displaymath}
Mimicking this, we can introduce a similar concept for a pre-conciliation 
$G$ in $X$, called a {\it treaty} at $x$, by 
\begin{displaymath}
 G_x \eqdef \bigcup_{W_j\ni x} G(W_j)/\mathord{\approx_{\rho}}
 = \varinjlim_{W_j\ni x}G(W_j)
\end{displaymath}
where the limit is the inductive limit induced by mediations:
\begin{displaymath}
 \begin{array}{l}
  \displaystyle
  G(W_0) \stackrel{\delta^{U_0}_{U_1}}{\longrightarrow} G(W_1)
   \stackrel{\delta^{U_1}_{U_2}}{\longrightarrow} G(W_2)
   \stackrel{\delta^{U_2}_{U_3}}{\longrightarrow} 
   \cdots
   \varinjlim_{W_j\ni x} G(W_j)\\[2mm]
  W_0\subseteq W_1\subseteq W_2\subseteq \ldots
 \end{array}
\end{displaymath}
Here, the equivalence relation $\approx_{\delta}$ is defined by the following equivalence: for every $W_1, W_2\in\Fer(X), t_i\in W_i$
\begin{displaymath}
 t_1\approx_{\delta} t_2
 \;\iff\; 
 \exists W\supseteq W_1\cup W_2: \delta^W_{W_1}(t_1) = \delta^W_{W_2}(t_2).
\end{displaymath}
Elements of a {\it treaty} $G_x$ at $x$ are equivalence classes
(the analog of what we call {\it germs} in sheaf theory).
We propose to call such a class an {\it agreement} at $x$.
Then classes of {\it agreements} constitute {\it treaties} and
a pre-conciliation can be constructed from such treaties. 

Starting from any pre-conciliation, we can build an associated 
{\it conciliation}. 
Let $G$ be a {\it pre-conciliation}, we can consider the set $G^a$ 
of the following maps defined on each closed set: 
\begin{displaymath}
 a_W : W \rightarrow \bigcup_{x\in X}G_x: x \mapsto a_W(x)
\end{displaymath}
such that $a_W(x)\in G_x$. 
We are showing that $G^a$ is in fact a conciliation.
The {\it agreement} $a(x)$ is an equivalence class of the {\it treaty}
$G_x$, namely an element $[(s,W)]$, where
$s\in G(W)$, $x\in W$ and
$[(s,W)] = [(t,V)] \iff s \approx_{\delta} t$, if $x\in V$. 
One can then define, in a very natural way, the mediations of $G^a$
by the following formula: 
\begin{displaymath}
 \ddelta{G^a}{W}{W_1}(a_W(x))
 =
 \ddelta{G^a}{W_2}{W_1}\bigl([(s,W_1)]\bigr)
 \eqdef \bigl[\bigl(\ddelta{G}{W_2}{W_1}(s),W_2\bigr)\bigr]
\end{displaymath}
with $W_1\subseteq W_2$.
One can then prove that $G^a$ is a {\it conciliation}.
This process, generating a conciliation from a preconciliation, could serve as the starting point to define the quotient of two conciliations $G$ and 
$H$.
Giving the map $W\rightarrow G(W)/H(W)$ whenever it has a meaning,
we get a preconciliation $Q$.
One can then define $G/H$ as the associated conciliation $Q^a$.

\section{Conciliations and paraconsistent logic of closed sets}

First of all let us introduce the concept of Brouwer algebra
(see \cite{McKinsey1946}) which is the dual concept of Heyting algebra (see \cite[pp.~129--134]{Awodey2010}).
A Brouwer algebra, or co-Heyting algebra, denoted by
$\struct{B,\wedge, \vee,0,1,\psdiff}$, is a bounded distributive lattice with minimal element $0$ and maximal element $1$, endowed with a binary operation $\psdiff$, called “pseudo-difference”, satisfying the following property:
\begin{displaymath}
 \forall a,b,c\in B: (a \psdiff b)\leq c \iff a \leq(b \vee c)
\end{displaymath}
where $\leq$ is the order relation defined on the lattice by:
$a\leq b \iff a=a\wedge b \iff b = a\vee b$.

We know that Boolean algebra is related to classical propositional logic.
We also know that Heyting algebra is related to intuitionistic
propositional calculus.
Brouwer algebra corresponds to paraconsistent propositional logic 
(see \cite{McKinsey1946}, \cite{Rasiowa1974}, \cite{Goodman1981}).
The correspondence goes as follows.
To any proposition $p$ of this paraconsistent logic one can associate an
element $a(p)$ belonging to the Brouwer algebra
$\struct{B,\wedge, \vee,0,1,\psdiff}$.
To the conjunction and to the disjunction, one can associate the binary 
operations $\wedge$ and $\vee$, respectively. 
The contradiction $\bot$ and the tautology $\top$ are associated to the
elements $0$ and $1$ respectively.
Finally, the paraconsistent negation $\neg$ corresponds to 
$1 \psdiff a(p)$ in the Brouwer algebra.
We recall that paraconsistent logic does not satisfy the non-contradiction
principle, i.e., $p\wedge \neg p \leftrightarrow \bot$ .   

A classical example of Brouwer algebra is built from the set of closed sets of a topological space.
We can check that 
$(\Fer(X),\cap,\cup,\empty,X,\Psdiff)$ is a Brouwer algebra,
where the pseudo-difference is defined as follows:
\begin{displaymath}
 (Z\Psdiff W) \eqdef \overline{\complement_XW}\cap Z
\end{displaymath}
(where $\complement_XW$ denotes the complement of $W$ in $X$ and 
where the bar denotes the closure of a set, i.e. the least closed set containing this set).
In this algebraic interpretation, the paraconsistent negation $\neg p$
corresponds to
\begin{displaymath}
 (X \Psdiff W(p)) \eqdef \overline{\complement_X W(p)} \cap X
  = \overline{\complement_XW(p)}.
\end{displaymath}

As the space $X$ let us take the real line and let us define $\Fer(X)$
as the set of its closed intervals.
We immediately check that the noncontradiction principle 
$p\wedge \neg p \leftrightarrow \bot$ is not satisfied.
Indeed, if we choose $W(p)=[0,1]$ we get: 
\begin{displaymath}
 W(p) \cap \overline{\complement_XW(p)} = \{0\} \cup \{1\} \neq \emptyset.
\end{displaymath}

Intuitively we can say that contradiction is confined only on the border of the closed set.
This confinement prevents contradictions to invade all logical deductions and gives to this paraconsistent logic a non-trivial character.   

In fact one can prove that any Brouwer algebra is a subalgebra of a Brouwer
algebra of closed sets of a particular topological space
(see \cite{McKinsey1946}).
This result corresponds to another one that is well-known in the theory of Heyting algebras saying that every Heyting algebra is a subalgebra of a Heyting algebra of open sets of a topological space.
We know also this classical result stating that every Boolean algebra can be viewed as a subalgebra of a Boolean algebra of subsets of a set. 

Let us consider the set
$\Cov(W_i) = \{Z\in \Fer(X) \st Z \supseteq W_i\}$
defined above.
Following James and Mortensen (see \cite{Mortensen1995}), let us show that 
$\struct{\Cov(W_i),\cap,\cup,\emptyset,X,\PSdiff{W_i}}$ is a Brouwer 
algebra, where the pseudo-difference is given by: 
\begin{displaymath}
 (Z\PSdiff{W_i} W) \eqdef (\overline{\complement_XW} \cap Z) \cup W_i
\end{displaymath} 
for all $Z,W \in \Cov(W_i)$.

We have to prove that, for every $ Z,W,T \in \Cov(W_i)$,
\begin{displaymath}
 (Z\PSdiff{W_i} W) \subseteq T
 \iff
 Z \subseteq (W\cup T),
\end{displaymath}
i.e.:
\begin{displaymath}
 (\overline{\complement_XW}\cap Z) \cup W_i \subseteq T
 \iff 
 Z\subseteq (W\cup T).
\end{displaymath}
But this is easily checked by set considerations. 
Therefore $\Cov$ is a conciliation in $X$.
Let us emphasize the fact that in this definition we take the union with
$W_i$ in order that $(Z\PSdiff{W_i} W)$ belongs to $\Cov(W_i)$.
In fact, in general, if $Z,W\in \Cov(W_i)$ we do {\it not} have:
\begin{displaymath}
 (\overline{\complement_XW} \cap Z) = (Z\Psdiff W)\in \Cov(W_i).
\end{displaymath} 

The conciliation $\Cov$ in $X$ is a map which associates to any closed 
set $W$ of $X$, a particular Brouwer algebra built with all closed sets
containing $W$.
The definition of the conciliation $\Cov$ is coherent with the structure of 
the Brouwer algebra due to the fact that mediations are preserving the 
pseudo-difference $\PSdiff{W_i}$, 
i.e.
\begin{displaymath}
 \forall W_i\subseteq W_j:\;
 \delta^{W_j}_{W_i}(Z\PSdiff{W_i} W)
 \;=\;
 \delta^{W_j}_{W_i}(Z)\PSdiff{W_j} \delta^{W_j}_{W_i}(W).
\end{displaymath}
The proof of this statement follows from a result previously obtained by James et Mortensen and described in an unpublished book (\cite{James}):
\begin{displaymath}
\begin{array}{l}
 \displaystyle \delta^{W_j}_{W_i}(Z \PSdiff{W_i} W)
   = \bigl( (\overline{\complement_XW}\cap Z ) \cup W_i\bigr)\cup W_j
   = (\overline{\complement_XW}\cap Z ) \cup W_j
   = Z \PSdiff{W_i} W \\[2mm]
 \displaystyle \delta^{W_j}_{W_i}(Z) \PSdiff{W_j} 
   \delta^{W_j}_{W_i}(W)
   = \bigl( \overline{\complement_X(W\cup W_j)} \cap (Z\cup W_j)\bigr)
     \cup W_j
\end{array}
\end{displaymath}

Then, we check that:
\begin{align*}
   \bigl(\overline{\complement_X(W\cup W_j} \cap (Z\cup W_j)\bigr)
      \cup W_j 
  & = \overline{\bigl(\complement_X(W\cup W_j) \cap (Z\cup W_j)\bigr)
      \cup W_j}  \\
  & = \overline{\bigl(\complement_X(W\cup W_j)\cup W_j\bigr)
         \cap \bigl((Z\cup W_j)\cup W_j\bigr)}\\
  & = \overline{\bigl(\complement_X(W\cup W_j)\cup W_j\bigr)
         \cap (Z\cup W_j)}\\
  & = \overline{\bigl((\complement_X(W)\cap \complement_X(W_j))\cup W_j\bigr)
         \cap (Z\cup W_j)}\\
  & = \overline{\bigl((\complement_X(W)\cup W_j)
  	     \cap (\complement_X(W_j)\cup W_j)\bigr)
         \cap (Z\cup W_j)}\\
  & = \overline{(\complement_X(W)\cup W_j)
         \cap (Z\cup W_j)}\\
  & = \overline{(\complement_X(W)\cap (Z\cup W_j))
         \cup (W_j\cap (Z\cup W_j))}\\
  & = \overline{(\complement_X(W)\cap Z)\cup 
         (\complement_X(W) \cup W_j) \cup W_j}\\
  & = \overline{(\complement_X(W)\cap Z)\cup W_j}\\
  & = \overline{(\complement_X(W)\cap Z)}\cup W_j\\
  & = Z \PSdiff{W_j} W.
\end{align*}

\section{The category of conciliations}

If we want to compare conciliations in $X$, it is natural to introduce a notion of morphism of conciliations.
Let $G$ and $H$ be conciliations on $X$.
A morphism of conciliations $G\rightarrow H$ can be defined using the following commutative diagram:
\begin{displaymath}
 \begin{CD}
  G(Z) @>>> H(Z) \\
  @A{\textstyle\ddelta{G}{Z}{W}}AA  @AA{\textstyle\ddelta{H}{Z}{W}}A\\ 
  G(W) @>>> H(W)
 \end{CD}
\end{displaymath}

This means that the morphism is compatible with mediations of both conciliations (we have distinguished mediations of both conciliations with the help of an index at the left of the mediation symbol).  

The set of conciliations in $X$ and the morphisms of conciliations generate objects and arrows of a category: the category of conciliations on 
$X$ denoted by $\Ccs{X}$.
The axioms defining a category can be easily verified.
For any conciliation, there exists an identity morphism and the composition of two morphisms is still a morphism.
Furthermore, this composition is associative. 
The conciliation $1:W_i\rightarrow \{*\}\restr{W_i}$, 
together with the trivial mediations
$\ddelta{1}{W_j}{W_i}(*\restr{W_i}) = *\restr{W_i}$,
happens to be a terminal object of $\Ccs{X}$.
For any conciliation $G$, it exists effectively a unique morphism
$G\rightarrow 1$ compatible with all mediations, i.e. rendering the following diagram commutative:
\begin{displaymath}
 \begin{CD}
  G(W_j) @>>> 1(W_j) \\
  @A{\textstyle\ddelta{G}{W_j}{W_i}}AA  
  @AA{\textstyle\ddelta{1}{W_j}{W_i}}A\\ 
  G(W_i) @>>> 1(W_i)
 \end{CD}
\end{displaymath}

It is interesting to wonder what are the subobjects of $\Ccs{X}$.
What does a subconciliation mean? 
A subconciliation $G'$ of $G$ is a conciliation satisfying the following commutative diagram, for all $W_i\subseteq W_j$:
\begin{displaymath}
 \begin{array}{ccc}
   G'(W_j)   &  \subseteq & G(W_j)   \\[1mm]
   \vcenter{\llap{$\ddelta{G'}{W_j}{W_i}$}}\Big\uparrow   & & 
       \Big\uparrow\vcenter{\rlap{$\ddelta{G}{W_j}{W_i}$}}\\[3mm]
   G'(W_i)   &  \subseteq & G(W_i)
 \end{array}
\end{displaymath}

This diagram means that inclusions remain always compatible with mediations. 
We are showing that $\Ccs{X}$ can be characterized by the following commutative diagram:
\begin{displaymath}
 \begin{array}{ccc}
   G'(W_i)   &  \stackrel{\mskip7mu j}{\subseteq} & G(W_i)   \\[1mm]
   \Big\downarrow   & & 
     \Big\downarrow\vcenter{\rlap{$\psi_{G'}$}}\\[3mm]
   1(W_i) & \underset{\phi}{\longrightarrow} & \Cov(W_i)  
 \end{array}
\end{displaymath}
where $G'$ is a subconciliation of $G$.
The maps $\psi_{G'}$ and $\phi$ are defined as follows: 
\begin{itemize}
\item 
  $\psi_{G'}(j(t_i)) = Z_{t_i}
  = \bigcap \big\{Z_k \st Z_k\supseteq W_i 
  \AND t_i\in G'(W_i)
  \AND \ddelta{G}{Z_k}{W_i}(j(t_i))\in G'(Z_k)\big\}$.

  The image of this map is thus the least closed set containing $W_i$ 
  such that the mediation $\ddelta{G}{Z_k}{W_i}$ of $G$ lifts
  $j(t_i)$ inside $G’(W_i)$. 

\item $\phi:1(W_i)\rightarrow \Cov(W_i): * \mapsto W_i$.

  This maps associates to $*$ the minimal element of $\Cov(W_i)$.
  When $G'$ is a subconciliation of $G$, one checks that:
  $\forall t_i\in G'(W_i): \psi_{G'}(j(t_i))=\phi(*)=W_i$. 
\end{itemize}

For two conciliations $G_1$ et $G_2$, one can define their product, 
associating to each closed set $W_i$ the Cartesian product of 
$G_1(W_i)$ and $G_2(W_i)$.
Given these conciliations, one can define the set of morphisms relating 
both.
The preceding commutative diagram shows also that $\Ccs{X}$ possess
an object $\Cov$ classifying subconciliations.
This category has thus the properties of a topos. 
As $\Cov(W_i)$ can be endowed with a structure of Brouwer algebra (and not 
of Heyting algebra as in an ordinary topos), 
$\Ccs{X}$ is an example of {\it complement topos}
(see Mortensen and Lavers in \cite[p.~105]{Mortensen1995})
(one could call it a co-Heyting or Brouwer topos); namely a topos whose classifying object is a Brouwer algebra.
Let us note that in fact $\Ccs{X}$ gives an example of co-Heyting 
topos which is not a bi-Heyting topos (see \cite{Reyes1996}).

\section{Towards a cohomology of conciliations}

The aim of this section is to show that conciliation generates cohomology in the same way that sheaf theory generates sheaf cohomology
(see \cite{Perrin1995}, \cite{Neeman2007}, \cite{Ward1990}).
Let $\{W_j\}_{j\in I}$ be a closed covering of $X$ and let $G$ be a 
conciliation (of rings) in $X$, which associates to any closed set
$W_j$ a ring $G(W_j) = \struct{A,+,\cdot}$. 
Now we define a {\it $p$-cochain with values in this conciliation,}
as the map that associates to any union of closed sets
$W_{j_0}\cup W_{j_1}\cup\ldots\cup W_{j_p}$,
$j_{\alpha}\in I$, $j_0< j_1< \ldots < j_p$, an element
$t_{j_0,j_1,\ldots,j_p}$  of
$G(W_{j_0}\cup W_{j_1}\cup\ldots\cup W_{j_p})$.
A $p$-cochain is thus a set 
$\{t_{j_0,j_1,\ldots,j_p}\}_{j_0< j_1< \ldots < j_p}$.
The set of $p$-cochains will be denoted by $C^p(I,G)$.
A $0$-cochain is a set $\{t_j\}_{j\in I}$, $t_j\in G(W_j)$.
A $1$-cochain is a set $\{t_{j_0,j_1}\}_{j_0<j_1}$,
$t_{j_0,j_1}\in G(W_{j_0}\cup W_{j_1})$. Etc.

One can define a {\it coboundary operator} $D$ such that:
\begin{displaymath}
 D: C^p(I,G) \rightarrow C^{p+1}(I,G):
 \bigl\{t_{j_0,j_1,\ldots,j_p}\bigr\}
 \mapsto
 \bigl\{(Dt)_{j_0,j_1,\ldots,j_p,j_{p+1}}\bigr\}
\end{displaymath}
where
\begin{displaymath}
 (Dt)_{j_0,j_1,\ldots,j_p,j_{p+1}}
  = \sum^{p+1}_{k=0} (-1)^k
    \delta^{W_{j_0}\cup W_{j_1}\cup \ldots \cup 
    W_{j_{p+1}}}_{W_{j_0}\cup W_{j_1}\cup \ldots \cup \widehat{W_{j_k}}
    \cup \ldots \cup W_{j_{p+1}}}
    t_{j_0,j_1,\ldots,\widehat{j_k},\ldots,j_p,j_{p+1}}
\end{displaymath}
($\widehat{\kern 1ex}$ means that you have to cancel the symbol).
We check that $D$ satisfies the property of a good coboundary operator, namely: $DD=0$.  

If we are computing the action of the coboundary operator on the
$0$-cochains, we get: 
\begin{displaymath}
(Dt)_{j_0,j_1}
    = \delta^{W_{j_0}\cup W_{j_1}}_{W_{j_1}}\: t_{j_1}
      - \delta^{W_{j_0}\cup W_{j_1}}_{W_{j_0}}\: t_{j_0}
\end{displaymath}
and if we put this to zero, we are lead to the condition (\ref{prop-vi})
satisfied by all elements of $G(W_j)$ if $G$ is a conciliation.
Similarly, we can compute the action of the coboundary operator on the
$1$-cochains:
\begin{displaymath}
 (Dt)_{j_0,j_1,j_2}
 = \delta^{W_{j_0}\cup W_{j_1}\cup W_{j_2}}_{W_{j_1}\cup W_{j_2}}
     \: t_{j_1,j_2}
   - \delta^{W_{j_0}\cup W_{j_1}\cup W_{j_2}}_{W_{j_0}\cup W_{j_2}}
     \: t_{j_0,j_2}
   + \delta^{W_{j_0}\cup W_{j_1}\cup W_{j_2}}_{W_{j_0}\cup W_{j_1}}
     \: t_{j_0,j_1}.
\end{displaymath}

Now, let us suppose that the $1$-cochain is defined by $t=Ds$.
The preceding formula becomes:
\begin{displaymath}
 \begin{array}{rcl}
  (Dt)_{j_0,j_1,j_2}
    & = & \delta^{W_{j_0}\cup W_{j_1}\cup W_{j_2}}_{W_{j_1}\cup W_{j_2}}
          \Big(\delta^{W_{j_1}\cup W_{j_2}}_{W_{j_2}} \: s_{j_2}
          - \delta^{W_{j_1}\cup W_{j_2}}_{W_{j_1}} \: s_{j_1} \Big)\\[4mm]
    & - &\delta^{W_{j_0}\cup W_{j_1}\cup W_{j_2}}_{W_{j_0}\cup W_{j_2}}
          \Big( \delta^{W_{j_0}\cup W_{j_2}}_{W_{j_2}} \: s_{j_2}
          - \delta^{W_{j_0}\cup W_{j_2}}_{W_{j_0}} \: s_{j_0} \Big)\\[4mm]
    & + & \delta^{W_{j_0}\cup W_{j_1}\cup W_{j_2}}_{W_{j_0}\cup W_{j_1}}
          \Big( \delta^{W_{j_0}\cup W_{j_1}}_{W_{j_1}} \: s_{j_1}
          - \delta^{W_{j_0}\cup W_{j_1}}_{W_{j_0}} \: s_{j_0}\Big)
 \end{array}
\end{displaymath}
Distributing médiations and applying the property (\ref{prop-iv}) we get:
\begin{displaymath}
 \begin{array}{rcl}
 (Dt)_{j_0,j_1,j_2}
 & = & \delta^{W_{j_0}\cup W_{j_1}\cup W_{j_2}}_{W_{j_2}}
       \: s_{j_2}
        - \delta^{W_{j_0}\cup W_{j_1}\cup W_{j_2}}_{W_{j_1}}
          s_{j_1}\\[4mm]
 & - & \delta^{W_{j_0}\cup W_{j_1}\cup W_{j_2}}_{W_{j_2}}
          \: s_{j_2}
       + \delta^{W_{j_0}\cup W_{j_1}\cup W_{j_2}}_{W_{j_0}}
         s_{j_0}\\[4mm]
 & + & \delta^{W_{j_0}\cup W_{j_1}\cup W_{j_2}}_{W_{j_1}}
         \: s_{j_1}
       - \delta^{W_{j_0}\cup W_{j_1}\cup W_{j_2}}_{W_{j_0}}
         \: s_{j_0}
 \end{array}
\end{displaymath}
this allows us to check that $DDs$ is zero.
As in the case of sheaves, we can introduce cohomology groups of
conciliation setting:
$H^p(I,G) = Z^p(I,G)/B^p(I,G)$ 
where $Z^p(I,G)$ is the set of $p$-cocycles,
i.e. of $p$-cochains $t$ such that
$Dt=0$ and where $B^p(I,G)$ is the set of $p$-coboundaries,
i.e. of $p$-cochains $b$ such that there exists $(p-1)$-cochain $s$
such that $b=Ds$.
The definition of the cohomology group $H^p(I,G)$ is motivated
by the fact that the $p$-cocycle $t$ is only given up to a 
coboundary $b=Ds$, because $Dt=0=D(t+b)=D(t+Ds)$.

\section{Globales and universales}

We want to suggest now that the classical notions of {\it locales} and 
{\it quantales} have some dual analogs in the framework of conciliations.

A {\it locale} is a complete lattice L satisfying the following property: 
\begin{displaymath}
   x \wedge \bigvee Y = \bigvee(x\wedge Y),
   \qquad x\wedge Y \eqdef \{x\wedge y \st y \in Y\AND Y\subseteq L\}.
\end{displaymath}

In other words, $L$ is a {\it “meet infinitely distributive”} complete distributive lattice. 

A {\it quantale}
$\struct{Q,\wedge,\vee,\bullet,0,1}$  is a structure defined on a set $Q$
such that
\begin{enumerate}
\item $\struct{Q,\wedge,\vee}$ is a complete lattice with nul element $0$    
   and universal
   element $1$,

\item $\struct{Q,\bullet}$ is a monoid with identity $1'$ for the operation
   $\bullet$ 

\item the operation $\bullet$ is distributive with respect to $\vee$
   in the following sense:
   \begin{displaymath}
     \forall a\in Q\, \forall S\subseteq Q:
         a\bullet \bigvee S = \bigvee\{a\bullet s \st s \in S\}
         \AND
         \bigl(\bigvee S\bigr)\bullet a = \bigvee\{s\bullet a \st s \in S\}
   \end{displaymath}
   where $\bigvee S = \{\bigvee_j s_j\st s_j\in S\}$.
\end{enumerate}

One can immediately check that a locale is a quantale if we take 
$\forall a\in Q: a\bullet b \eqdef a\wedge b$. 

One can prove that a {\it quantale} such that the operation $\bullet$ 
is idempotent ($\forall a\in Q: a\bullet a = 1'$) and such that the 
universal element coincide with the monoid identity ($1=1’$) is a 
{\it locale}.

Let us now introduce the concept of {\it globale} as a complete lattice
$L$ satisfying the following property:
\begin{displaymath}
 w\vee \bigwedge Z = \bigwedge(w \vee Z)
\end{displaymath}
where $w\vee Z \eqdef \{w\vee z \st z\in Z \AND Z\subseteq L\}$.

In other words, $L$ is a {\it “join infinitely distributive”} complete
 lattice. 

An example of {\it globale} is given by the lattice of closed sets,
$\Fer(X)$, of a topological space $X$.
Now, we define 
$\forall w_1,w_2\in\Fer(X), w_1\wedge w_2 \eqdef w_1\cap w_2$.
As the union of an arbitrary number of closed sets is not necessarily 
closed, we have to define  
$\forall w_1\in\Fer(X): \bigvee_j w_j\eqdef \overline{\bigcup_j w_j}$
using the closure of the set union. 
Both operations endowed $\Fer(X)$ with a structure of complete lattice. 
We give here a proof of the well-known property saying that this lattice is effectively a {\it “join infinitely distributive”} lattice:

First of all we have:
\begin{displaymath}
 \forall w\in\Fer(X), z\in Z\subseteq\Fer(X):
 w\cup z \supseteq w \cup \bigcap Z.
\end{displaymath}
We can deduce that: 
\begin{displaymath}
 \bigwedge(w\vee Z) = \bigcap (w\cup Z)
 \supseteq w\cup \bigcap Z = w\vee \bigwedge Z.
\end{displaymath}

Let us note that we can write:
$w\cup \bigcap Z = w\vee \bigwedge Z$
because $\bigcap Z$ is a closed set and the union of two closed sets
is still a closed set. 

We have now to prove that:
$w\vee \bigwedge Z \supseteq \bigwedge(w\vee Z)$.
 
We have in fact
$\forall z\in Z: w\cup z \supseteq \bigcap(w\cup Z)$.
Defining $u=\bigcap(w\cup Z)$ we find 
$w\cup z\supseteq u$.
But we have also: 
\begin{displaymath}
 \begin{array}{l}
   z = (w\cup z)\cap \complement_W(w-z)
   \supseteq u \cap \complement_X(w-z)
   = u \cap \complement_X(w\cap \complement_Xz)
   = u \cap (\complement_X\cup z)\\
     \Rightarrow\; z\supseteq (u\cap \complement_Xw)
     \cup (u\cap z).
 \end{array}
\end{displaymath}

Taking the intersection $v=\bigcap Z$ of elements $z\in Z$, we get:
\begin{displaymath}
 \begin{array}{l}
  v\supseteq (u\cap v)\cup(u\cap \complement_Xw)
  = [(u\cap v)\cup u]\cap[(u\cap v)\cup \complement_Xw]
  = u\cap (u\cup \complement_Xw)\cap(v\cup \complement_Xw)\\
  \Rightarrow v\supseteq u\cap (v\cup\complement_Xw)\\
  \Rightarrow v\cup w \supseteq [u\cap (v\cup\complement_Xw)]\cup w
    = (u\cup w)\cap(v\cup\complement_Xw\cup w)
    = u \cup w \supseteq u\\
  \Leftrightarrow w\cup v \supseteq u
  \Leftrightarrow w\cup\bigcap Z \supseteq \bigcap(w\cup Z).   
 \end{array}
\end{displaymath}

The lattice of closed sets is thus a globale. 

Next we define a kind of “dual” version of the quantale.
We will call it a {\it universale}.
As a matter of definition one could propose the following statement:
$\struct{P,\wedge,\vee,0,1,\circ}$ is a {\it universale}
if the set $P$ is such that: 
\begin{enumerate}
 \item $\struct{P,\wedge,\vee,0,1,\circ}$ is a complete lattice with
   zero element $0$ and universal element $1$.

 \item $\struct{P,\circ}$ is a monoid with identity element $0'$ with   
   respect to the 
   operation $\circ$.
 
\item $\forall b\in P, \forall T\subseteq P:
   b\circ \bigwedge T = \bigwedge(b\circ T)
   = \bigwedge\{b\circ t\st t\in T\};
   (\bigwedge T)\circ b = \bigwedge(T\circ b)$.
\end{enumerate}

One can show that such a structure has the following properties:

The operation $\circ$ is such that
$\forall a,b,c\in P: a\leq b \Rightarrow a\circ c \leq b\circ c$,
where the order relation is the relation generated by the lattice:
$a\leq b\iff a\wedge b; b = a\vee b$.
We have indeed,
\begin{displaymath}
 a\leq b \iff a = a\wedge b
 \Rightarrow a\circ c = (a\wedge b)\circ c
 = (a\circ c)\wedge (b\circ c) 
 \iff a\circ c \leq b\circ c
\end{displaymath}

If the zero element is identified with the identity element of the monoid
operation ($0=0’$), we get: 
$\forall a,b\in P: a \circ b \geq a\vee b$.
This is easily proved:
\begin{displaymath}
 \begin{array}{l}
  a\geq 0 = 0' \Rightarrow
  a\circ b \geq 0'\circ b = b
  \iff a\circ b \geq b \iff b = (a\circ b)\wedge b\\
  b\geq 0 = 0' \Rightarrow
  a\circ b \geq a\circ 0' = a
  \iff a\circ b \geq a \iff a = (a\circ b)\wedge a
 \end{array}
\end{displaymath}
So
\begin{displaymath}
 \begin{array}{l}
  a\vee b = [(a\circ b)\wedge a] \vee [(a\circ b)\wedge b]\\
  \iff a\vee b = (a\circ b) \wedge (a\vee b)\\
  \iff a\vee b \leq a\circ b.
 \end{array}
\end{displaymath}

We can exhibit the link between the {\it globale} and the {\it universale}
by the following property: 
if the operation of the {\it universale}
$\struct{P,\wedge,\vee,0,1,\circ}$ is idempotent and if the zero element of
the lattice is identified with the identity element of the monoid,
then $\struct{P,\wedge,\vee,0,1}$ is a {\it globale}. 

If $\struct{P,\wedge,\vee,0,1}$ is a {\it globale} and if one identifies
$\circ$ to $\vee$, one gets 
$\forall a\in P: a\circ a = a \vee a = a$ and this shows the idempotency
property.
Furthermore, one check that
\begin{displaymath}
 \forall a\in P: a\vee 0 = a =a \circ 0' = a \vee 0' \Rightarrow 0=0'
\end{displaymath}
and that
\begin{displaymath}
 a\circ \bigwedge_j b_j = a \vee \bigwedge_j b_j =
 \bigwedge_j(a\vee b_j) = \bigwedge_j(a \circ b_j). 
\end{displaymath}
This shows that $\struct{P,\wedge,\vee,0,1}$ is a {\it universale} with an
idempotent monoid operation and with the zero element of the lattice 
identified with the identity element of the monoid. 

If now $\struct{P,\wedge,\vee,0,1,\circ}$ is a {\it universale} with 
idempotent monoid operation and with the universal element of the lattice
identified with the identity element of the monoid, then if we identify 
$\circ$ to $\vee$ we get: 
\begin{displaymath}
 \forall a \in P, \forall T\subseteq P:
 a\circ \bigwedge T = \bigwedge(a\circ T)
 \iff a\vee \bigwedge T = \bigwedge(a\vee T)
\end{displaymath}
and this shows that $\struct{P,\wedge,\vee,0,1}$ is a {\it globale}.

\section{Speculations}

The main motivation of this paper was to build a formal model of (social, political, legal) situations where one reaches a conciliation of contradictory positions. 
This leads us to propose axioms defining a structure on a topological space that we called {\it conciliation}.
It was surprising to realize that such a structure was in fact a kind of dualization of a sheaf;
the duality here (which has not to be taken in the strict sense of category theory) has to be understood as the reversion of the role of global and local levels combined with the change from open to closed sets.
In the case of a sheaf, we glue “locally” some objects defined on open
sets, checking a coherence constraint on the intersections of open sets and 
then, we can reconstruct a well-defined “global” (i.e. at the level of the union of open sets) object of the same type.
On the contrary, in the case of a conciliation, we conciliate globally some objects defined on closed sets, checking a coherence constraint on the union of closed sets, and then we can reconstruct well-defined “local” objects 
(i.e.~at the level of intersections of closed sets) of the same type. 

In other words, starting from trivial local objects, a sheaf allows us to 
build non-trivial global object by a gluing procedure. At the opposite, 
starting from trivial global objects, a conciliation allows us to generate 
some non-trivial local object.     

Our paper may also help to shed some new light on a very important diagram occurring in the book of Shahn Majid, {\it Foundations of Quantum Group Theory} \cite{Majid1995}, but also in some of his papers (see for example 
\cite{Majid2008}) and which illustrates the progress towards a search of a 
quantum gravity theory regulated by the meta-principle of self-duality expressing his representation-theoretic philosophy. 

In the upper part of this diagram (see the figure below), we find concepts
and theories related to quantum physics.
We find, for example, Heyting algebra, which is the starting point for an 
understanding of quantum logic in which one drops the law of excluded middle.
We find also quantales.

In the lower part of the diagram, one can find notions that are connected with general relativity and gravitation.
Majid locates, in this part, co-Heyting algebras, i.e.~Brouwer algebras. 
This is because the latter are related to logics in which one can drop the non-contradiction principle.
As noted by Lawvere (see \cite{Lawvere1991}), this fact enables us to 
introduce an operator
$\partial p \eqdef p\wedge \neg p \neq 0$  which behaves as a derivation operator 
and which satisfies Leibniz rule:
$\partial(p\wedge q) = (\partial p \wedge q) \vee (p\wedge \partial q)$.
As derivation is the main ingredient to build vector fields in differential 
geometry, one can think that Brouwer algebra has something to do with the 
birth of geometry and step after step with general relativity and 
gravitation.

The upper and lower part of the diagram are related by duality.
For example, if we replace open sets by closed sets and the intuitionistic 
implication
\begin{displaymath}
 (U\Rightarrow V) = (\overdot{(\complement_XU)} \cup V)
\end{displaymath}
($\overdot{\complement_XU}$ being the interior of
$\complement_XU$) by the paraconsistent pseudo-difference
$(Z\Psdiff T)=(\overline{\complement_XT}\cap Z)$, we pass from the 
Heyting algebra in the upper part to Brouwer algebra in the lower part.
The horizontal line separating both parts of the diagram corresponds to self-dual structures.
An example of such structure is the Boolean algebra which can be considered as an algebra that is at the same time a Heyting and a Brouwer algebra.
The set of subsets of a set is a Boolean algebra and its elements are at 
the same time open and closed sets, what the logicians called clopen sets.
This is once again a facet of duality.

It is possible to improve Majid’s inspiring diagram adding on the 
horizontal line “bi-Heyting algebras”
(see \cite[pp.~249, 300--301]{Awodey2010}) which are at the same time 
Heyting and Brouwer algebras without being equivalent to a Boolean algebras.
One example of such algebras is the lattice of subgraphs of a given graph 
(see \cite[pp.~29--31]{Reyes1996}).
It is also possible to add in the upper part {\it locales} which are 
complete Heyting algebras and which constitute generalizations of the 
concept of topological space (see \cite[pp.~472--475]{MacLane1992}).
One example of {\it locale} is the lattice of open sets of a topological 
space (see \cite{Isham2005}).
These {\it locales} can pave the way to the {\it quantales}
(see \cite{Mulvey1986}, \cite{Coniglio2000}) quoted explicitly by Majid,
and which are generalizations of {\it locales} well suited to formalize 
quantum theories. 
An example of such {\it quantales} is given by the lattice of ideals of a 
commutative ring with unity.
The upper part of Majid’s diagram is the natural place for the mathematics 
of sheaves.
We can for instance build sheaves on a locale and the category of such 
sheaves is called a {\it localic topos} (see \cite[pp.~472-527]{MacLane1992}), which generalizes the topos of sheaf on a usual topological space.

Using the “duality à la Majid”, we can now come back to the lower part of the diagram.
The dual of a {\it locale} is in fact a complete Brouwer algebra and this is 
what we have called a {\it globale}.
A natural example of such a structure is the lattice of closed sets of a 
topological space.
In fact properties of such lattices have been studied many years ago, in particular in the {\it Séminaire Dubreil} (see \cite{Lesieur1953}).
We can also consider the dual of a {\it quantale} and we find what we have 
called above a {\it universale}.
The dual notion of a {\it sheaf on a locale} is now a {\it conciliation on
a globale} and it would be interesting to study the category of such 
conciliations: what we might call a {\it globalic topos}.
This topos is maybe a complement topos but this remains to be proved.
It is also possible to think about the category of conciliations on a 
{\it globale} but up to now we don’t know exactly what it is.
Using the Majid’s duality philosophy we could try to find self-dual 
structures that are at the same time {\it locale} and {\it globale} or at 
the same time {\it quantale} and {\it universale}.
But we don’t know if this is really relevant.
We will try to clarify these problems in a forthcoming publication.

\setlength{\unitlength}{1mm}
\begin{center}
\begin{picture}(120,85)(0,-45)

\put(55,34){\makebox(0,0){\shortstack{Category of sheaves\\
on a locale: localic topos}}}
\put(55,24){\vector(0,1){5}}
\put(55,22){\makebox(0,0){Sheaf on a locale}}
\put(55,12.5){\vector(0,1){6}}

\put(20,11){\makebox(0,0){Quantale}}
\put(42,11){\vector(-1,0){8}}
\put(55,11){\makebox(0,0){Locale}}
\put(80,11){\vector(-1,0){8}}

\put(0,0){\makebox(0,0){?}}
\put(10,7){\vector(-4,-3){5}}
\put(14,0){\vector(-1,0){8}}
\put(10,-7){\vector(-4,3){5}}

\put(30,0){\makebox(0,0){?}}
\put(40,7){\vector(-4,-3){5}}
\put(44,0){\vector(-1,0){8}}
\put(40,-7){\vector(-4,3){5}}

\put(62,0){\makebox(0,0){bi-Heyting algebra}}
\put(87,7){\vector(-4,-3){5}}
\put(91,0){\vector(-1,0){8}}
\put(87,-7){\vector(-4,3){5}}

\put(110,11){\makebox(0,0){Heyting algebra}}
\put(110,2){\vector(0,1){6}}
\put(110,0){\makebox(0,0){Boolean algebra}}
\put(110,-2.5){\vector(0,-1){4}}
\put(110,-11){\makebox(0,0){\shortstack{Brouwer algebra\\
                                       (co-Heyting algebra)}}}

\put(20,-11){\makebox(0,0){\shortstack{Universale\\(co-quantale)}}}
\put(42,-11){\vector(-1,0){8}}
\put(55,-11){\makebox(0,0){Globale}}
\put(80,-11){\vector(-1,0){8}}

\put(55,-13.5){\vector(0,-1){5}}
\put(55,-23){\makebox(0,0){\shortstack{Conciliation (co-sheaf)\\
                                       on a globale}}}
\put(55,-28.5){\vector(0,-1){5}}
\put(55,-38){\makebox(0,0){\shortstack{Category of conciliations\\
                                       on a globale: globalic topos}}}

\end{picture}
\end{center}

\vspace{-2\baselineskip}

%
%
%
%
%

%
%
%
%
\begin{flushright}
 \begin{tabular}{l@{}}
  Dominique Lambert ({\tt d.lambert@fundp.ac.be})\\
  Bertrand Hespel ({\tt bertrand.hespel@fundp.ac.be})\\[3mm]
  D\'epartement Sciences-Philosophies-Soci\'et\'es\\
  Facult\'e des sciences\\
  Centre ESPHIN et NAXYS\\
  Universit\'e de Namur - FUNDP\\
  61, rue de Bruxelles -- B-5000 Namur (Belgium)
 \end{tabular}
\end{flushright}

\end{document}